\documentclass{article}
\usepackage{amsmath,amssymb,graphicx}
\usepackage[colorlinks,plainpages=false]{hyperref}

\hypersetup{colorlinks,%
citecolor=black,%
filecolor=black,%
linkcolor=black,%
urlcolor=black,%
}

\def\titl{Discrete analog of the Jacobi set for vector fields}
\def\auth{A.N. Adilkhanov, A.V. Pavlov, I.A. Taimanov}

\title{\titl\thanks{This work was supported by the Ministry of Education and Science of the Republic of Kazakhstan (program 0115PK03029) and Russian Foundation for Basic Research (grant 15-01-01671a).}
}

\author{
  A.N. Adilkhanov\footnote{
	National Laboratory ``Astana'', Nazarbayev University, 53, 
	Kabanbay Batyr Ave., Astana, 010000, Republic of Kazakhstan; 
	\href{mailto:aadilkhanov@nu.edu.kz}{aadilkhanov@nu.edu.kz}}, 
	A.V. Pavlov\footnote{
	North-Eastern Federal University, 677000 Yakutsk, Russia; 
	\href{mailto:av.pavlov@s-vfu.ru}{av.pavlov@s-vfu.ru}}, and 
	I.A. Taimanov\footnote{
	Sobolev Institute of Mathematics, 630090 Novosibirsk, 
	and Novosibirsk State University, 630090 Novosibirsk, Russia;
	\href{mailto:taimanov@math.nsc.ru}{taimanov@math.nsc.ru}}
}
\date{}

\hypersetup{pdftitle=\titl, pdfauthor=\auth, pdftoolbar=false,
plainpages=false, hyperindex=true, pdfdisplaydoctitle=true}

\begin{document}

\maketitle

\begin{abstract}
The Jacobi set is a useful descriptor of mutual behavior of functions defined on a common domain.
We introduce the piecewise linear Jacobi set for general vector fields on simplicial complexes. 
This definition generalizes the definition of the Jacobi set for gradients of functions
introduced by Edelsbrunner and Harer.


\end{abstract}

\section{Introduction}
\label{sec:1}

In this article we give a construction of a piecewise linear analog of the Jacobi set for vector
fields.  This set serves as a descriptor of the relation between vector fields defined on a common
domain.

For the gradient fields of Morse functions $f_1,\dots, f_k: D \to {\mathbb R}$, where $D$ is a
domain in ${\mathbb R}^N$, or more generally an $N$-dimensional manifold, the Jacobi set is the
subset of $D$  formed by all points at which the gradients of these functions are linearly
dependent.  This set  can be used for extracting useful information about the mutual behavior of
multiple functions \cite{measure}.  As Jacobi sets for a pair of functions on the plane it appears
for different reasons in~\cite{W1}  (see also \cite{W2}), and in general form it  was introduced in
\cite{EH}.

For applications, it is helpful to have a discrete analog of the Jacobi set, and such an analog for
functions defined on triangulated complexes was introduced in~\cite{EH}. In the same article, the 
problem of extending the proposed methods to general vector fields was posed.  We demonstrate how 
to do that on the example of pairs of vector fields on the plane.

\section{The piecewise linear Jacobi set}
\label{sec:2}

We recall the main definitions and results from \cite{EH}.

A \emph{Morse function} on a compact manifold $M$ is a function $f:M \to \mathbb{R}$ that has only
a~finite number of critical points, where the matrix of second derivatives is nondegenerate and the
function values are distinct from each other. For two Morse functions $f, g: M \to \mathbb{R}$
defined on a compact manifold $M$, their \emph{Jacobi set} is defined as the set $J(f, g)$ where
their gradients are linearly dependent. Equivalently, $J(f, g)$ can be described as the set of all
critical points of functions $h_\lambda = f+\lambda g$ and $e_\lambda = \lambda f + g$ for all
$\lambda \in \mathbb{R}$. For two generic Morse functions $f$, $g$ having no common critical points,
$J(f,g)$ is a 1-submanifold of $M$.

In the discrete case of two functions $f$, $g$ defined on the vertices of a~triangulation $K$ of a
$d$-dimensional manifold, we can extend $f$ and $g$ to PL-functions on the entire complex $K$ and
view each as a limit of a series of smooth functions. Motivated by this viewpoint, the discrete
Jacobi set ${\mathbb{J}(f,g)}$ is introduced~\cite{EH} as a 1-dimensional subcomplex of $K$
consisting of edges $uv$, with multiplicity, along which, in the limit, the critical points of
$h_\lambda$ and $e_\lambda$ travel as $\lambda$ varies.

To state the precise definition, we need some notation. Let $K$ be a simplicial complex. The
\emph{star} of its simplex $\sigma$ is the set of all simplices containing $\sigma$, and the
\emph{link} $\mathop{\mathrm{Lk}} \sigma$ consists of all simplices in the closure of the star of
$\sigma$ that are disjoint from~$\sigma$. Note that $\mathop{\mathrm{Lk}}\sigma$ is itself a
complex. Let $h$ be a real-valued function on the vertices of a simplicial complex $K$. For a
simplex $\sigma\in K$, define the lower link $\mathop{\underline{\mathrm{Lk}}} \sigma$ in $K$ with
respect to $h$ to be the portion of $\mathop{\mathrm{Lk}} v$ that bounds the set of all simplices in
the star of $v$ that have $v$ as the vertex with the maximal value of $h$.

Consider an edge $uv \in K$. We disregard the edges with $g(u)=g(v)$. Denote by
$\lambda^*=\lambda^*(u,v)$ the value of $\lambda$ that equalizes the values of the linear
combination $f+\lambda g$ at both ends of the edge: $f(u)+\lambda^* g(u) = f(v)+\lambda^* g(v)$.
Denote this linear combination by $h$: $h=f+\lambda^*g$. The link of the vertex $u$ is a
triangulation of a $(d-1)$-sphere containing~$v$. The multiplicity of the edge $uv$ is defined as
the sum of reduced Betti numbers $\tilde{\beta}_i$ of the lower link $\mathop{\underline{\mathrm{Lk}}} uv$ with
respect to $h$. The \emph{piecewise linear Jacobi set} of two functions $f$, $g$ on $K$ is defined
as the one-dimensional subcomplex $\mathbb{J}(f,g)$ of $K$ consisting of all edges having nonzero
multiplicity, together with their endpoints.

\begin{figure}[b]
\centering
\begin{tabular}{@{}c@{}c@{}}
\includegraphics{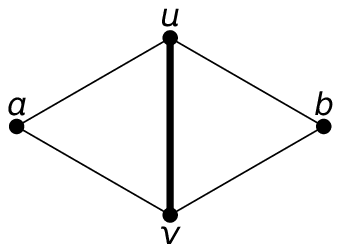} & \includegraphics{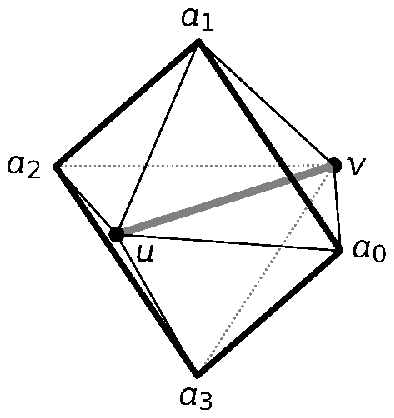}\\
(\textit{a})                & (\textit{b})
\end{tabular}
\caption{The link of an edge in a triangulation of (\emph{a}) 2-manifold (\emph{b})~3-manifold}
\label{figure:links}
\end{figure}

We now review the general definition for the special cases of 2- and 3-dimensional simplicial complexes.

In two dimensions, the star of an edge $uv$ consists of two 2-simplices neighboring along $uv$, and
$\mathop{\mathrm{Lk}}uv$ is just two vertices $a$, $b$ opposite $uv$ in these 2-simplices
(fig.~\ref{figure:links},~\emph{a}).  Thus, the edge $uv$ belongs to $\mathbb{J}(f,g)$ if and only
if the values $h(a)$, $h(b)$ are either both greater or both smaller than the value of $h(u)=h(v)$,
where $h=f+\lambda^*g$, and $\lambda^* = (f(u)-f(v))/(g(v)-g(u))$. This condition can be rewritten
in terms of function differences in adjacent vertices. For any function $f:K\to {\mathbb R}$ and
edge $xy$, denote $df(xy) = f(y)-f(x)$. Then we can write
\begin{equation} 
\begin{array}{rcl} 
\label{eq:JacobiEH}
uv \in \mathbb{J}(f,g) &
\Leftrightarrow & dh(ua) \text{ and } dh(ub) \text{ have the same sign,}\\ 
& & \text{where }h = f+\lambda^* g\text{~ and~} \lambda^* = df(uv)/df(vu).
\end{array} 
\end{equation}
This condition is actually symmetric in $u$ and $v$, since $dh(vx) = dh(ux) - dh(uv)$ for any vertex $x$, and $dh(uv)=0$.

In a triangulation of a three-dimensional manifold, the link $\mathop{\mathrm{Lk}}uv$ of an edge
$uv$ is a triangulation of a circle. The multiplicity of an edge $uv$ in the Jacobi set
$\mathbb{J}(f,g)$ is  equal to the sum of the reduced Betti numbers of the lower link
$\mathop{\mathrm{\underline{Lk}}}{uv}$  with respect to~$h$:
\[
\tilde{\beta}_{-1} + \tilde{\beta}_0+\tilde{\beta}_1,
\] 
where 
$\tilde{\beta}_{-1}$ is 1 if the lower link is empty, and 0 otherwise;
$\tilde{\beta}_0$ is one less than the number of connected components in
$\mathop{\mathrm{Lk}}\underline{uv}$ if this number is positive, and 0 otherwise;
$\tilde{\beta}_1$ is 1 if the lower link is the entire circle, and 0 otherwise.  

In fig.~\ref{figure:links},~\emph{b}, the link of $uv$ is shown in bold lines. This link is a
triangulation of a circle. Denote its consecutive vertices by  $a_0$, $a_1$, \ldots, $a_{k-1}$, and
put $a_k=a_0$. As previously, $dh(xy)$ stands for the difference $h(y)-h(x)$. Count the number of
times the difference $dh(ua_i)$  changes from negative to positive along the circle:
\[
 \beta_0 = \#\{i :~ dh(ua_i)<0 \text{~ and ~} dh(ua_{i+1})>0,~ 0\le i < k\}.
\]
Then the multiplicity of the edge $uv$ in $\mathbb{J}(f,g)$ is $|\beta_0-1|$. In particular, if the
lower link of the edge $uv$ with respect to $h$ consists of just a single component that is not the
entire circle, the edge does not belong to the Jacobi set.

\begin{figure}[h]
\centerline{\includegraphics[scale=.35]{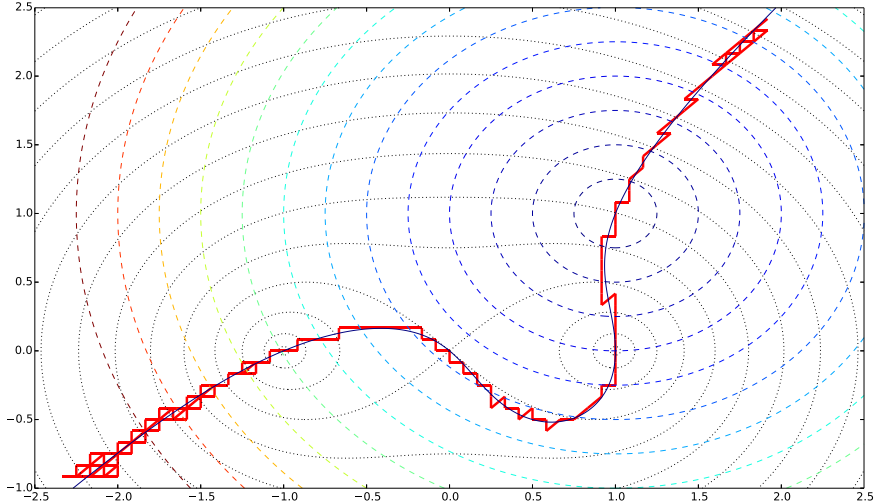}}
\caption{Zigzags in a simplicial Jacobi set (red) for functions $f = ((x-1)^2 + y^2)((x+1)^2 + y^2)$
(dotted level lines) and $g = (x-1)^2+(y-1)^2$ (dashed level lines). The triangulation is obtained 
from a square grid with step~$\frac{1}{6}$. The continuous black line is the smooth Jacobi set.}
\label{figure:zigzags}
\end{figure}

The functions $h_\lambda(v)=f+\lambda g$ are linear in $\lambda$ for any given $v$. Because of that,
any $v \in \mathop{\mathrm{Lk}} u$ changes its status as inside/outside the lower link of $u$ with
respect to $h_\lambda$  exactly once as $\lambda$ grows from $-\infty$ to $+\infty$, namely at
$\lambda=\lambda^*(u,v) = df(uv)/dg(vu)$. So in dimension 2, the number of connected components
$\beta_0$ in the lower link of $u$ is either the same at both extremes, or $0$ for one of them and
$1$ at the other. Obviously, for an edge $uv \notin \mathbb{J}$, passing $\lambda^*(u,v)$ does not
change $\beta_0$. For $uv \in \mathbb{J}$, passing $\lambda^*$ either changes $\beta_0$ by one or
does not change $\beta_0=1$ if on either side of $\lambda^*$ the lower link is all of the link
of~$u$. Counting the parity of $\beta_0$, we see that the number of edges $uv \in \mathbb{J}$ for a
fixed vertex~$u$, i.e. the degree of~$u$ in $\mathbb{J}$, must be even. A similar argument, after
unfolding each multiple critical point into multiple simple critical points, holds in any dimension:

\medskip
\textsc{Even-degree lemma} \cite{EH}. The degree of any vertex $u$ in $\mathbb{J}$ is even.
\medskip

Although the even-degree lemma guarantees that the discrete Jacobi set can be represented as a
continuous polyline, it may contain spurious cycles and zigzags, becoming incovienently large. For
example, if the simplicial complex in question is a fine enough regular triangulation of a plane,
the discrete Jacobi set may appear to fill entire 2-dimensional regions on the plane
(fig.~\ref{figure:zigzags}). A variety of simplification techniques exist for the smooth as well as
discrete versions of the Jacobi set (\cite{Natarajan}, \cite{Pascucci}).

\section{The piecewise linear Jacobi set for vector fields}
\label{sec:3}

The main idea behind our definition is as follows. The gradient 
\[
df  = \left(\frac{\partial f}{\partial x^1}, \dots, \frac{\partial f}{\partial x^N}\right)
\]
of a function $f: D \to {\mathbb R}$ is in fact a $1$-form which is a linear form on vector fields.
Indeed, its value for a vector field $X$ is the derivative of $f$ in the direction of $X$:
\[
D_X f = X^i \frac{\partial f}{\partial x^i},
\]
where we assume the summation over the repeated index. To obtain the gradient vector field we have
to raise the index by using some non-degenerate quadratic form $g^{ik}$ (usually the inverse of the
metric tensor $g_{ik}$):
\[
(\nabla f)^i = g^{ik}\frac{\partial f}{\partial x^k},
\]
where again we assume the summation over the repeated index $k$. The Euclidean metric is given by the tensor
\[
g_{ik} = g^{ik} = 
\begin{cases}
1 & \mbox{for $i=k$} \\
0 & \mbox{otherwise},
\end{cases}
\]
the gradient of the function and the gradient vector field look the same, but in general coordinates their
numerical expressions are different. We refer for detailed discussion, for instance, to \cite{NT}.

Since the lowering of the index (the convolution)
\[
X^i \to Y_k = g_{ik} X^i
\]
maps linearly dependent vector fields into linearly dependent $1$-forms, it is enough to define the Jacobi sets for $1$-forms.
  
For a triangulated complex $K$, $1$-forms $Y$ are linear functions on oriented $1$-chains, i.e., on oriented edges:
\[
Y(uv) \in {\mathbb R} \ \ \ \mbox{where $uv$ is an oriented edge in $K$}.
\]
We interpret the Edelsbrunner--Harer definition of the Jacobi set of two gradient vector fields as the definition
of the Jacobi set of two $1$-forms that are coboundaries of linear functions on the vertices of $K$:
\[
Y(uv)= df(uv) = f(v) - f(u).
\]
For a triangulation of a smooth manifold $K$ and a smooth function $f: K \to {\mathbb R}$ the
discretization of its gradient (covector) field is exactly given by the formula above where $f$ is
evaluated in the vertices of the triangulation.

Given a smooth $1$-form $\omega$ on a triangulated manifold, we have to construct a $1$-form on
oriented edges. The most natural way is to consider an edge as an oriented path and take an integral
of $\omega$ over the path:
\[
Y(uv) = \int_u^v \omega.
\]
For a smooth gradient field $df$ in Euclidean space we get
\[
Y(uv) = \int^u_v df = f(v)-f(u).
\]

A non-gradient vector field corresponds to a non-closed 1-form. Circular integrals of such a form
may not vanish, so generally, it is not true that $Y(uw) = Y(uv) + Y(vw)$.

Let $K$ be a simplicial complex that is a triangulation of a $d$-dimensional manifold, and $F$, $G$ be discrete 1-forms given by their values on all oriented edges $uv$ of $K$: 
\[
F(uv) = -F(vu),\quad  G(uv) = -G(vu).
\] 
Denote by $H_\lambda$ the linear combination $F+\lambda G$. For each edge $uv$ with $G(uv)\ne 0$, as previously, denote by $\lambda^*$ the value of the coefficient that makes this linear combination vanish along $uv$: 
\[
\lambda^* =\frac{F(uv)}{G(vu)}, \quad H_{\lambda^*}(uv)=0.
\]
For a vertex $w$ in $\mathop{\mathrm{Lk}}uv$ define the average of the values of the form $H_{\lambda^*}$ on the edges connecting $u$ and $v$ to $w$: 
\[
h(w) =\frac{1}{2}(H_{\lambda^*}(uw)+H_{\lambda^*}(vw)).
\]
Multiplicity of an edge $uv$ is defined as the sum of the reduced Betti numbers of the lower link of $uv$ with respect to $h$, and we define the \emph{Jacobi set $\mathbb{J}(F,G)$ of two discrete 1-forms $F$ and $G$} as the one-dimensional subcomplex
of $K$ consisting of all edges having nonzero multiplicity, together with their
endpoints.

\begin{figure}[b]
\centerline{\includegraphics[scale=1.1]{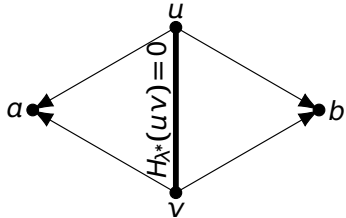}}
\caption{The edge test for 1-forms on a 2-dimensional complex}
\label{figure:uvab}
\end{figure}

In two dimensions ($d=2)$, this definition means that the Jacobi set of $F$ and $G$ consists of all edges $uv$ for which
the average of the values of $H_{\lambda^*}=F+\lambda^*G$ along $ua$ and $va$ has the same sign as the average of
its values along $ub$ and $vb$, where $a$ and $b$ are the two points of the link of $uv$ (fig.~\ref{figure:uvab}).
\begin{equation}\label{eq:JacobiV}
	uv \in \mathbb{J}(F,G) \Leftrightarrow (H_{\lambda^*}(ua)+H_{\lambda^*}(va))(H_{\lambda^*}(ub)+H_{\lambda^*}(vb))>0.
\end{equation}

Note that, as was the case for the condition (\ref{eq:JacobiEH}), this condition is also symmetric
in $u$, $v$. It is symmetric in $F$ and $G$ as well when all values of the forms $F$, $G$ on the 
edges are nonzero.

\begin{figure}[h!]
\centering 
\includegraphics[width=\textwidth]{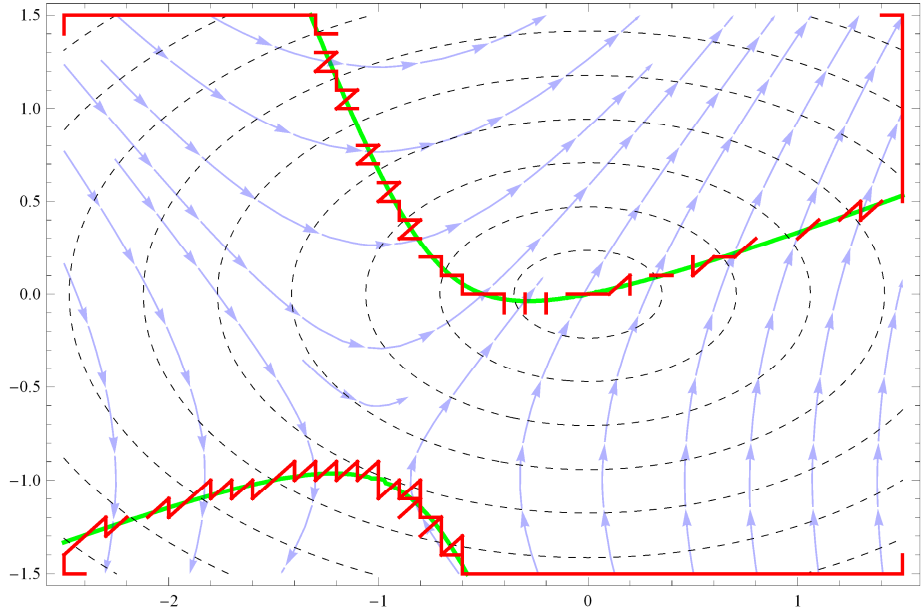}
\caption{Continuous Jacobi set (green line) and its piecewise linear approximation for 
 the forms $F(x,y)=(y+1)dx + 2(x+1)dy$ and $G(x,y)=(2x-3y)dx + (2x+3y)dy$, grid step size $h=0.1$}
\label{figure:sample1}
\end{figure}

However, the even degree lemma no longer holds for nongradient 1-forms. This is illustrated below
for the approximation of the Jacobi set for two smooth 1-forms on the plane (fig.~\ref{figure:sample1}). The
smooth Jacobi set is the set of points where the forms are linearly dependent, and is shown with
continuous green lines, while the piecewise linear Jacobi set for the triangulation of a square
grid with step size $h=0.1$ is shown in red.
 
Still, as can be seen in table~\ref{table:approx_var_he}, with the refinement of the grid the approximation converges to the smooth Jacobi set.

In applications, a vector field $X$  is usually given by its coordinates on a plane grid. A
reasonable approximation for the integrals of the corresponding 1-form $Y$ on the edges are scalar
products of the mean value of the vector field on the edge with the edge vector itself:
\[
Y(uv) = \frac{1}{2}(X_v+X_u, uv).
\]   

\begin{figure}[h!]
\centering	
\begin{tabular}{c}
\includegraphics{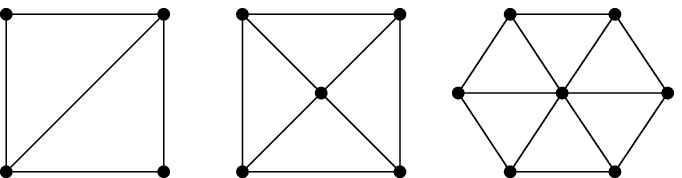} \\
(\emph{a}) T1 \hspace{44pt}
(\emph{b}) T2 \hspace{47pt}
(\emph{c}) T3 \hspace{10pt}\strut
\end{tabular}
\caption{Plane triangulations}
\label{figure:triangulations}
\end{figure}

We have also tested our definition for three different regular triangulations on the plane, shown in
fig.~\ref{figure:triangulations}: the diagonal grid $T_1$ (invariant with respect to rotations by
$\pi$), crossed $T_2$ (invariant with respect to rotations by $\pi/4$), and hexagonal $T_3$
(invariant with respect to rotations by $\pi/6$).

\begin{figure}
\noindent
\includegraphics[width=57mm]{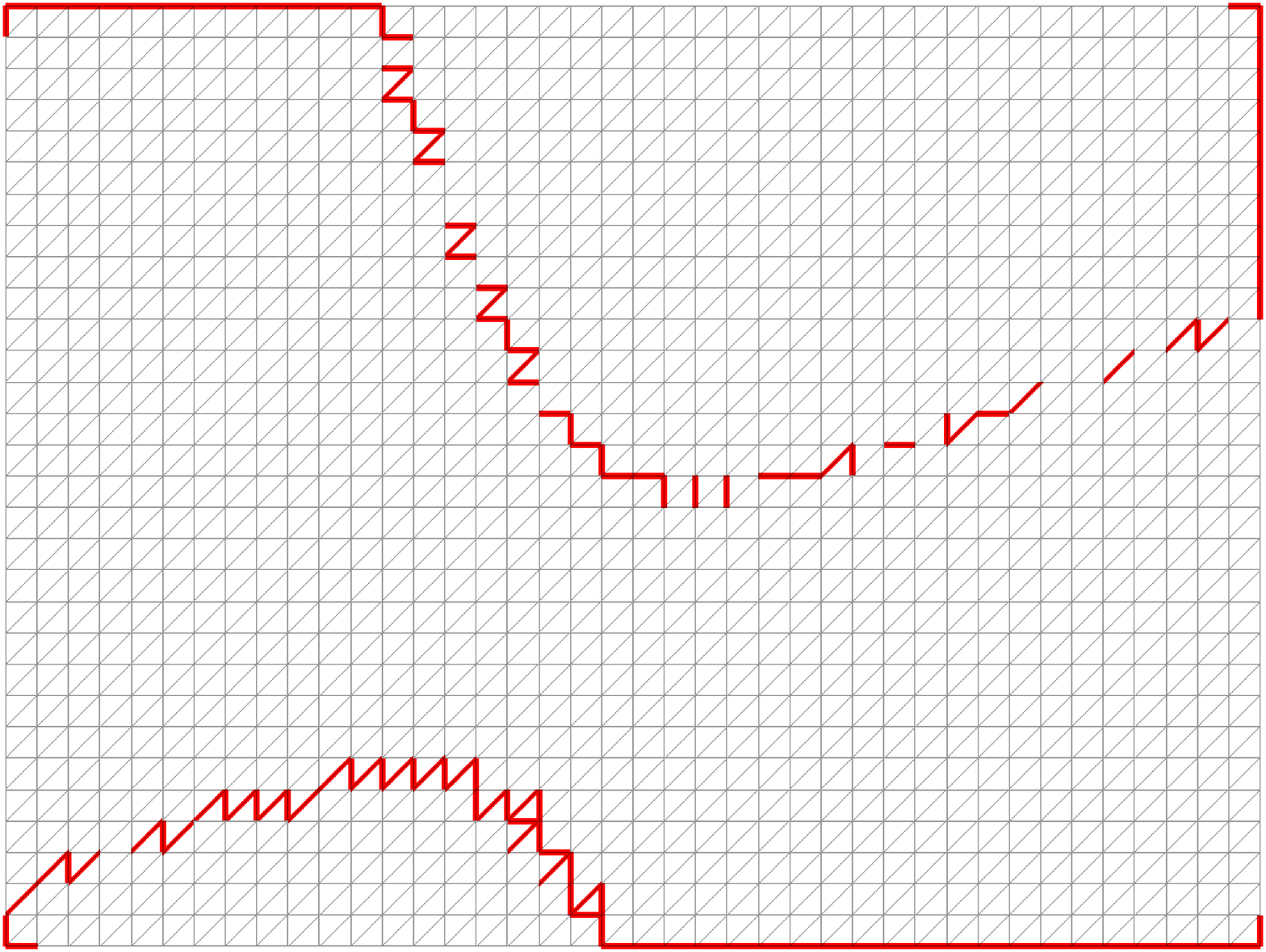}
\hfill
\includegraphics[width=57mm]{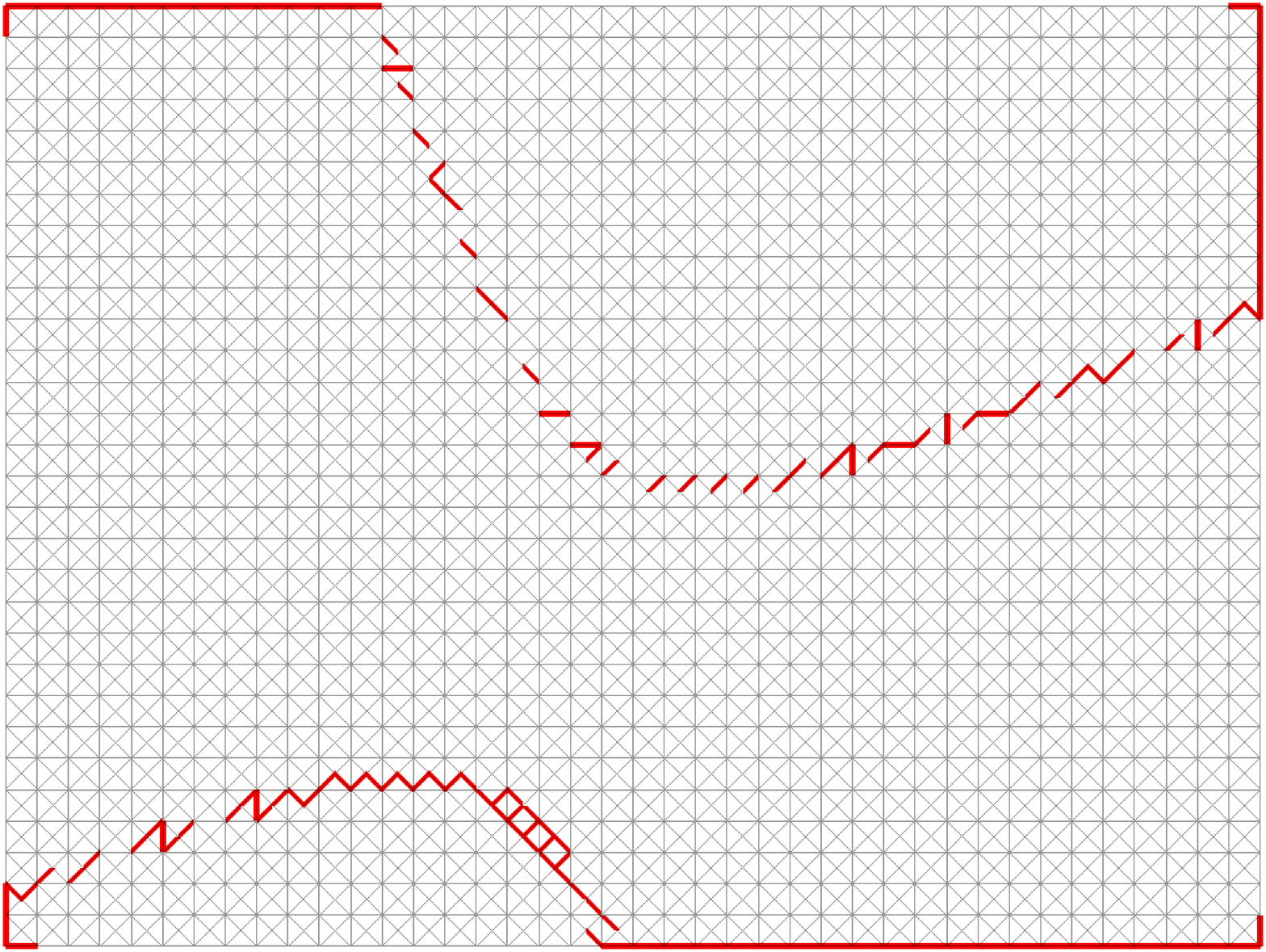}
\\[3mm]
\includegraphics[width=57mm]{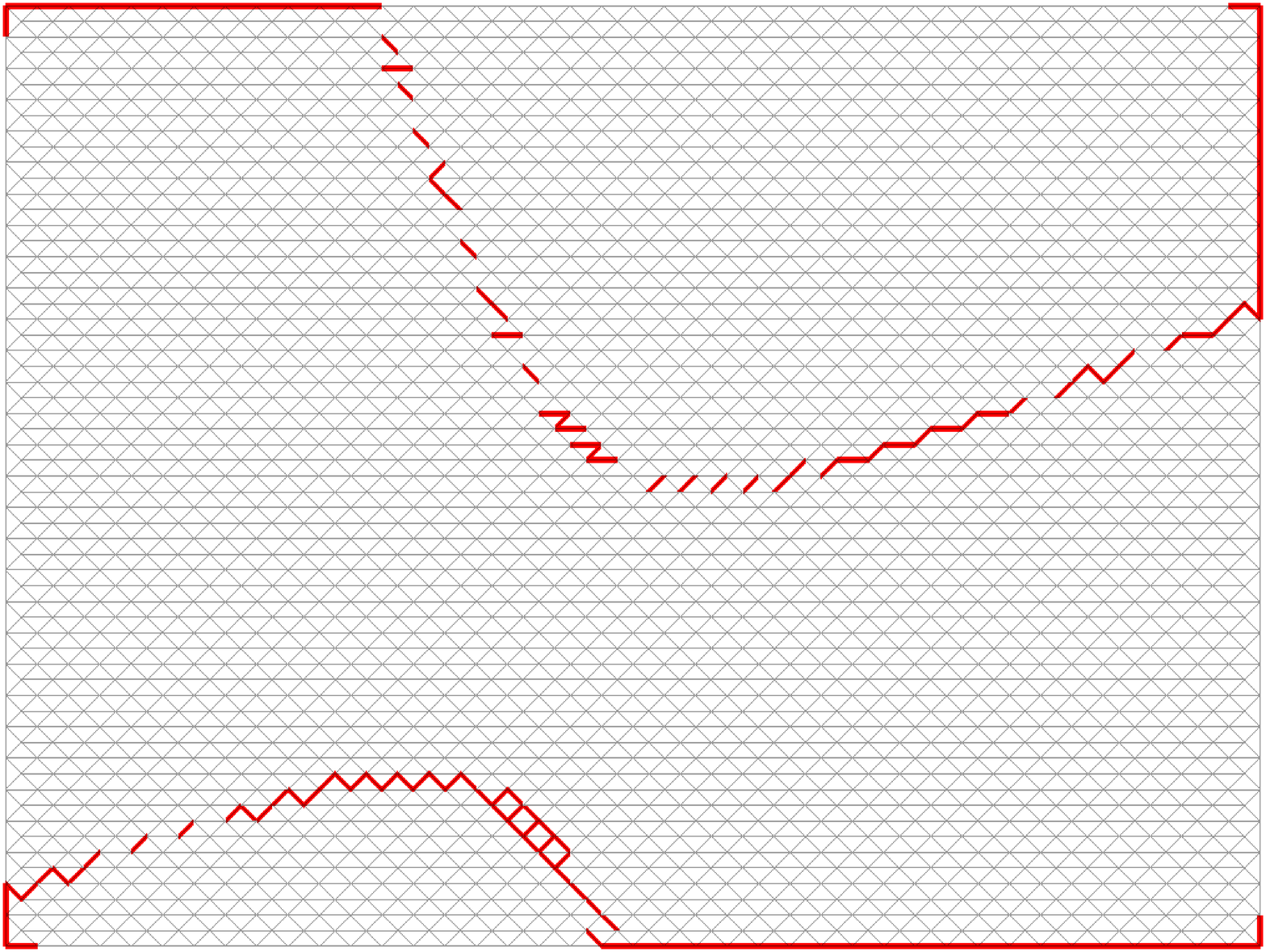}
 \caption{Piecewise linear Jacobi set of $F(x,y)=(y+1)dx + 2(x+1)dy$ and $G(x,y)=(2x-3y)dx + (2x+3y)dy$ for various triangulations on the plane}
 \label{figure:triangulationsJacobi}
\end{figure}

Results of these calculations are shown in fig.~\ref{figure:triangulationsJacobi}. As was the case for the Jacobi sets of Morse functions, the approximations differ, with no clear winner.

For better connectivity of the produced approximation, the edge test (\ref{eq:JacobiV}) can be modified to include cases where the absolute value of at least one of the factors is smaller than some threshold value $\varepsilon$:
\begin{equation}\label{eq:JacobiEps}
\begin{array}{r@{}l}
	 \mathbb{J}(F,G) = \{ & uv\in K \mid (H_{\lambda^*}(ua)+H_{\lambda^*}(va))(H_{\lambda^*}(ub)+H_{\lambda^*}(vb))>0.\\
	   & \text{or \:} |H_{\lambda^*}(ua)+H_{\lambda^*}(va)| < \varepsilon \;
	\text{\: or \:} |H_{\lambda^*}(ub)+H_{\lambda^*}(vb)| < \varepsilon \}
\end{array}
\end{equation}
This will improve the connectivity at the cost of thickening the Jacobi set. 

\begin{figure}
\noindent
\centerline{\includegraphics[width=.8 \textwidth]{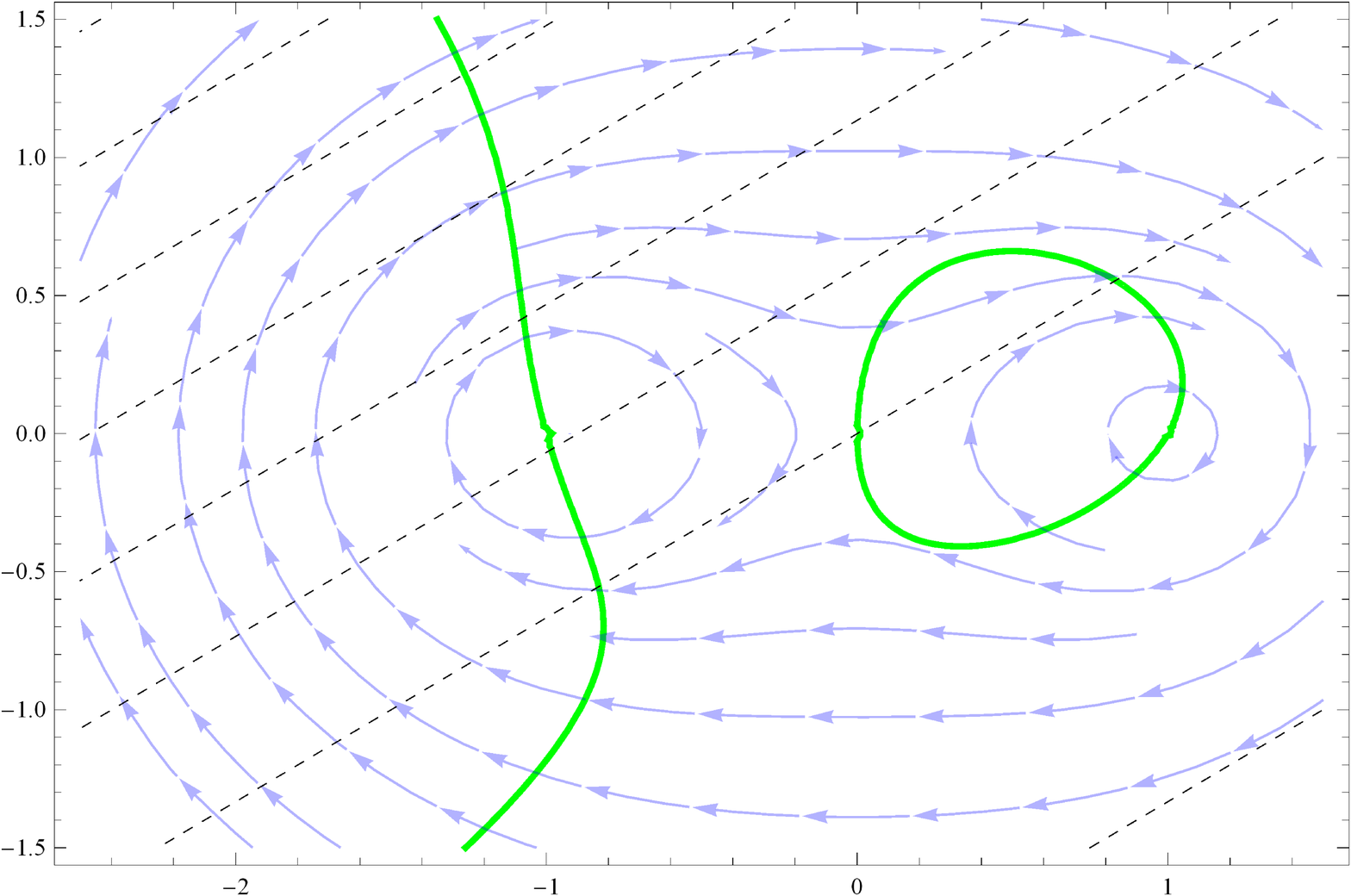}}
\caption{The fields $f$, $g$ and the smooth Jacobi set $J(F,G)$ for the corresponding forms\newline
$F(x,y) = y(x^2+y^2+1)dx - x(x^2+y^2-1)dy$,\quad $G(x,y) = (2x-3y-6)dx + (2x-3y)dy$ }
\label{fields}
\end{figure}

In fig.~\ref{fields}, we show the smooth Jacobi set, and in table~\ref{table:approx_var_he} illustrate the dependence of approximation, using the T1 triangulation scheme, on the grid step size $h$ and the threshold $\varepsilon$ in \eqref{eq:JacobiEps} for the forms 
\[
F(x,y) = y(x^2+y^2+1)dx - x(x^2+y^2-1)dy,~~ G(x,y) = (2x-3y-6)dx + (2x-3y)dy.
\]

\begin{figure}
\noindent
\centerline{\includegraphics[scale=.8]{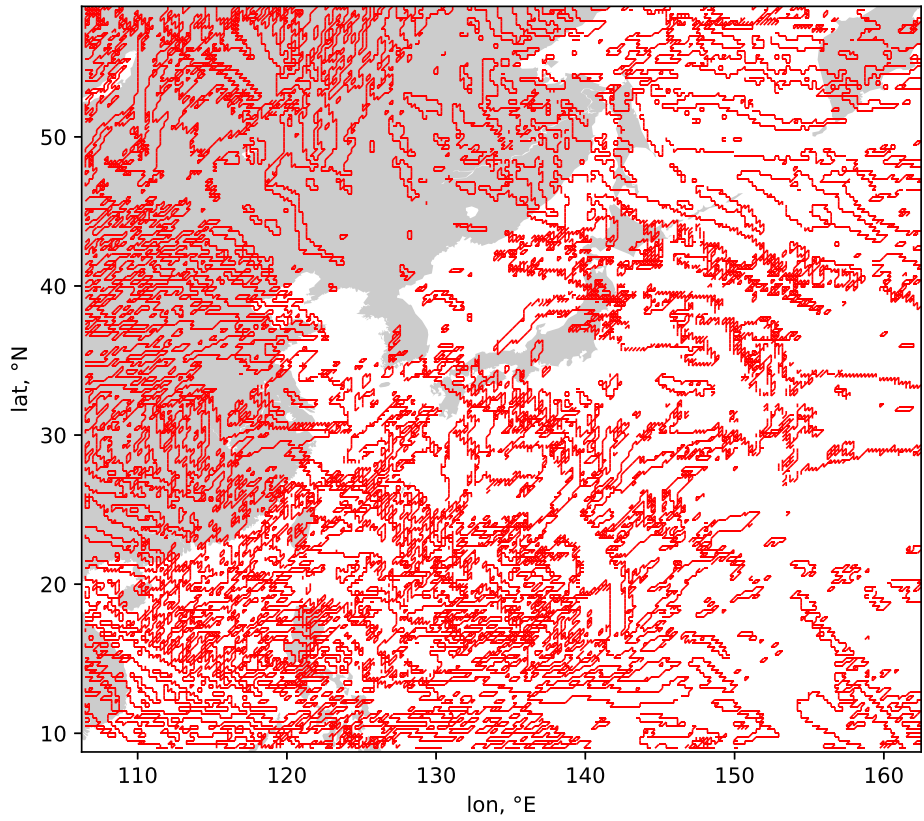}}
\caption{The discrete Jacobi set for the nongradient wind speed vector field
and the gradient of temperature over a part of Asia Pacific. GFS forecast data 
for 1 November 2018 used \cite{GFS}.}
\label{pacific}
\end{figure}

As in the case of the Jacobi set for functions, numerically approximated Jacobi set for vector fields may turn out to be very complicated. Sometimes, it might be an indication of a strong similarity between the vector fields, as in fig.~\ref{pacific}. However, it would be interesting to develop methods for its simplification similar to those proposed in \cite{Natarajan}, \cite{Pascucci}.

\begin{table}[p]
\label{table:approx_var_he}
\begin{center}
\begin{tabular}{|p{13mm}|c|c|c|}
\hline
 \centering$\varepsilon$ 
                 & $h=0.1$ & $h=0.05$ & $h=0.01$ \\
 \hline
 \raisebox{8mm}{$0.0$}         
               &\includegraphics[scale=0.08]{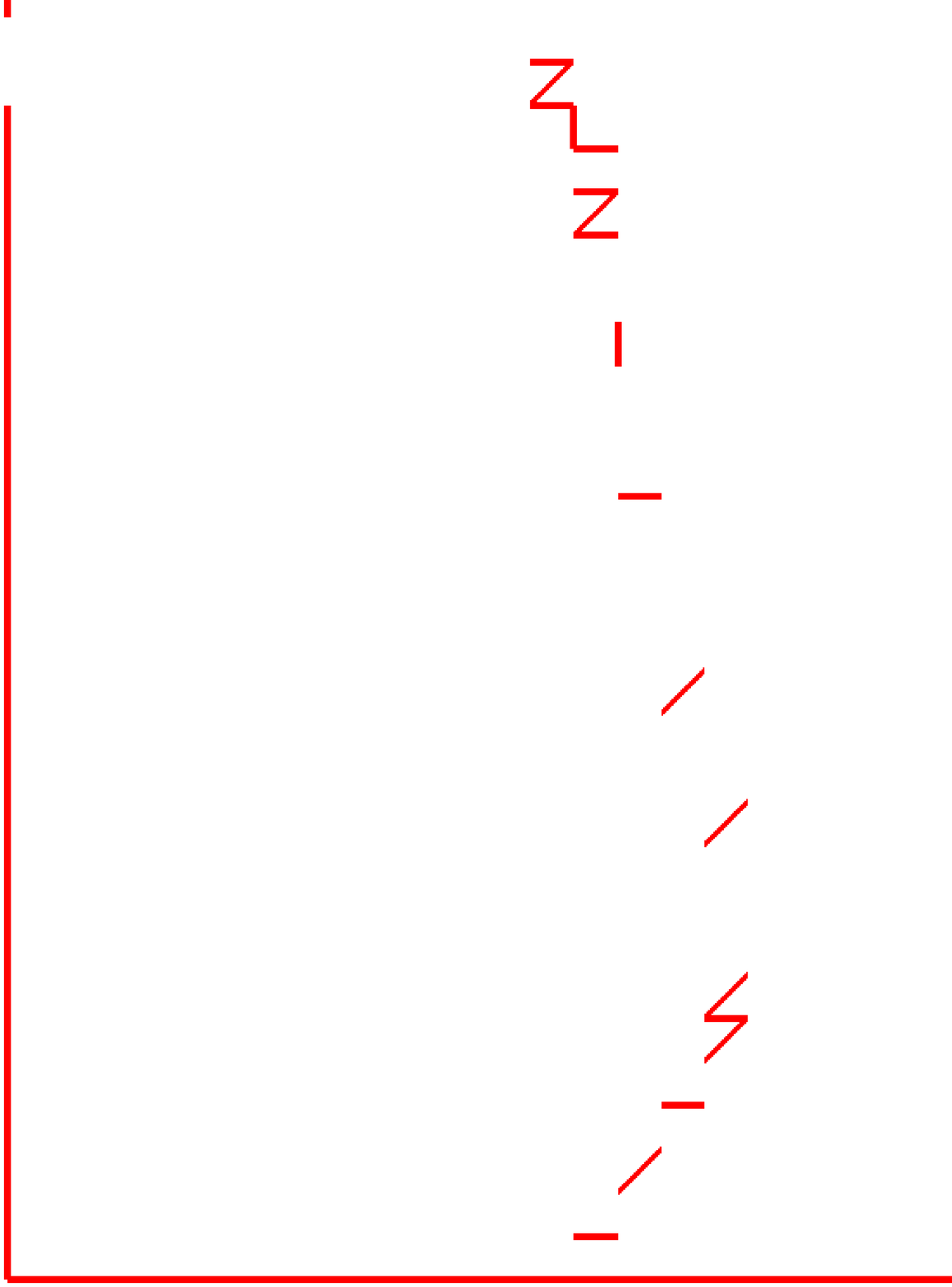}      
               &\includegraphics[scale=0.08]{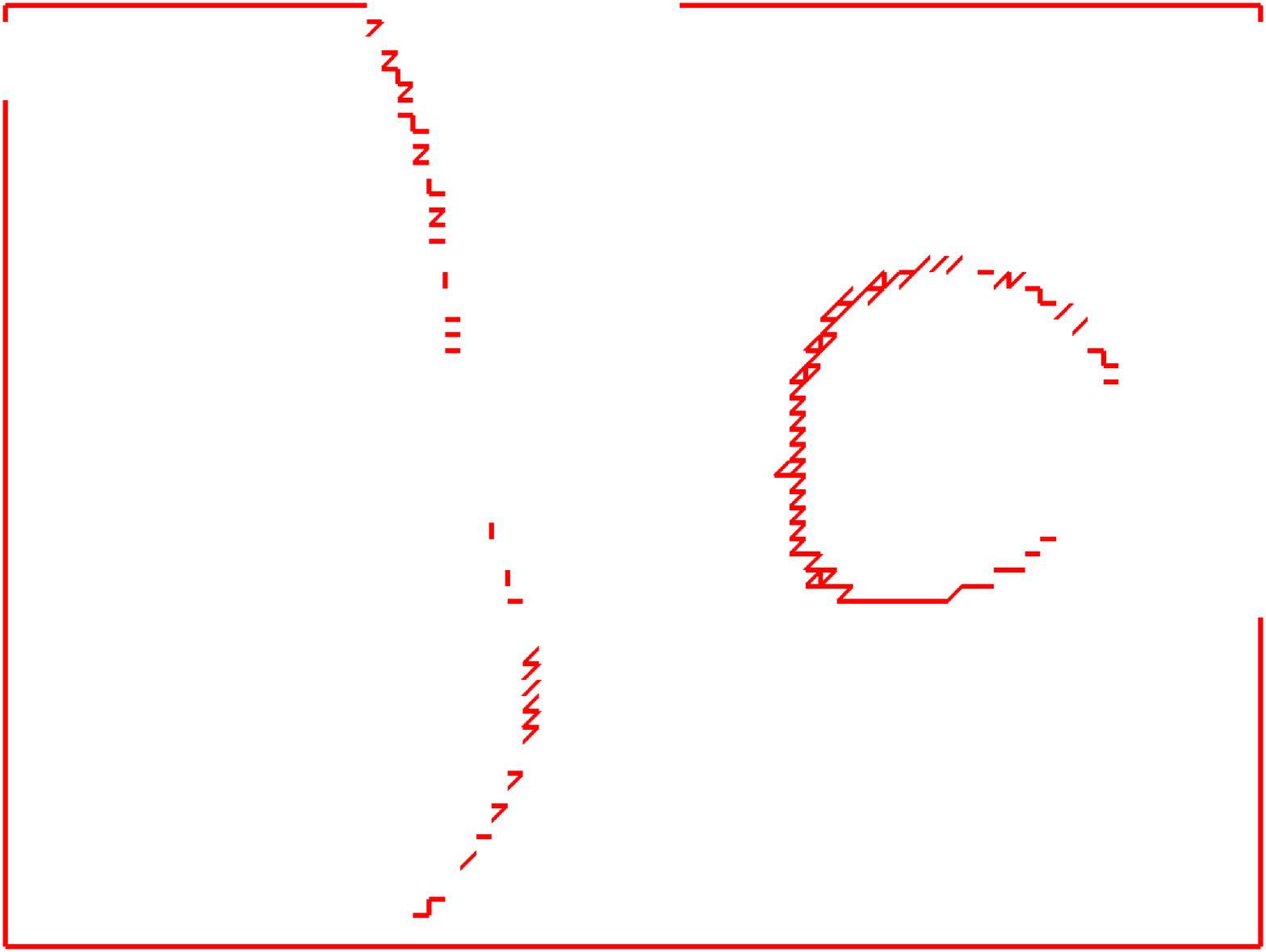}
               &\includegraphics[scale=0.08]{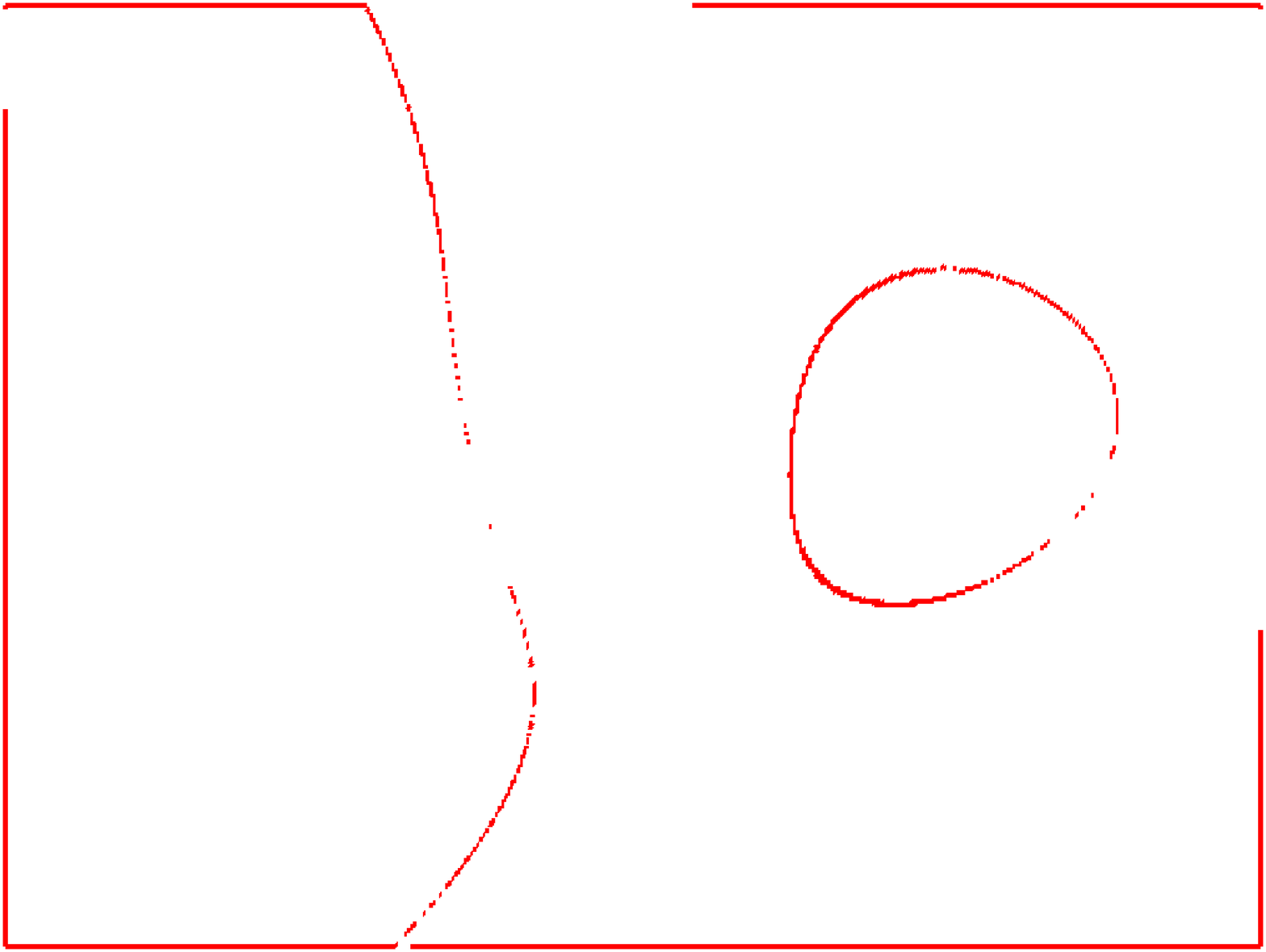}
               \\
 \raisebox{8mm}{$0.00005$}     
               &\includegraphics[scale=0.08]{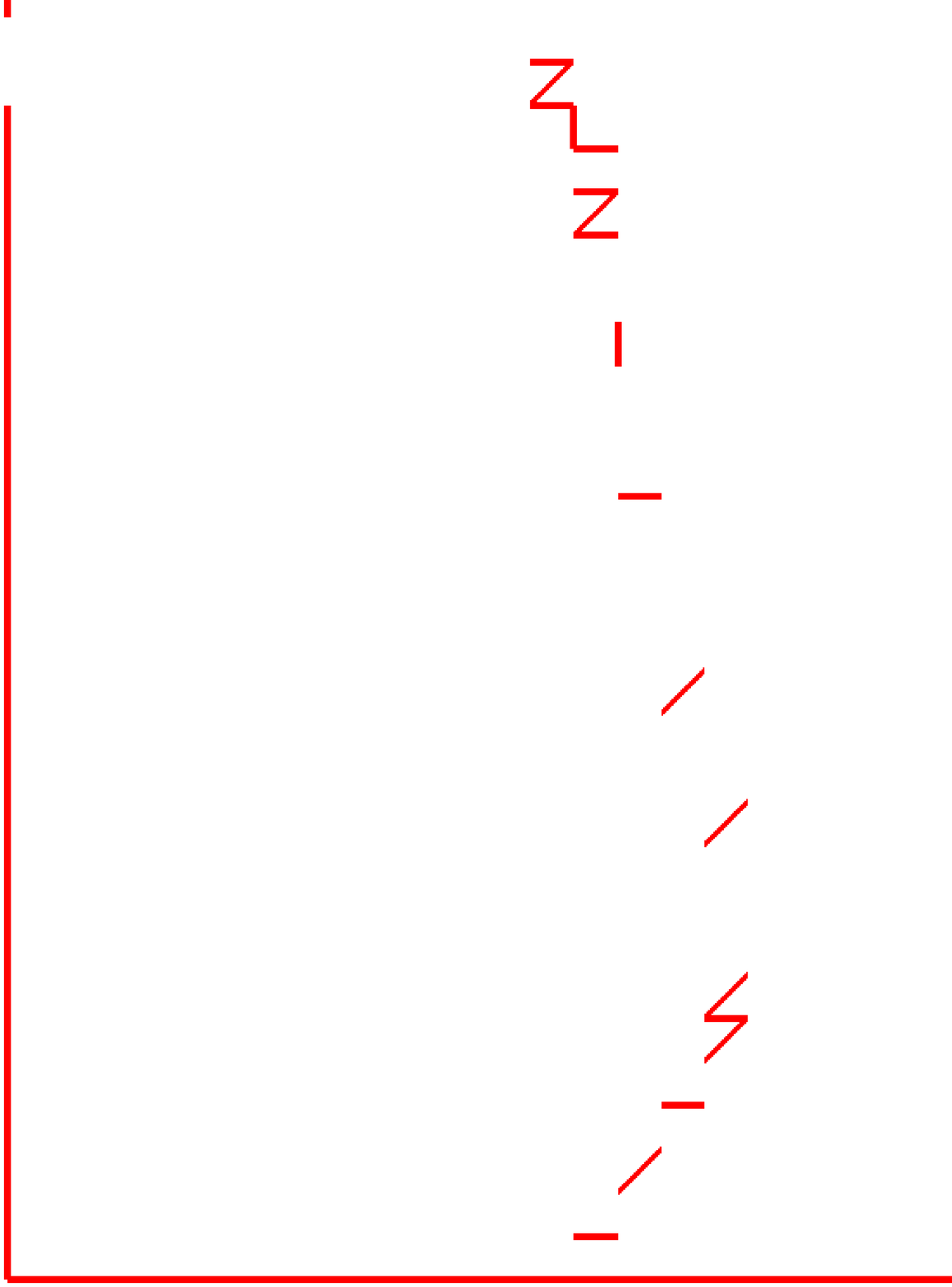}      
               &\includegraphics[scale=0.08]{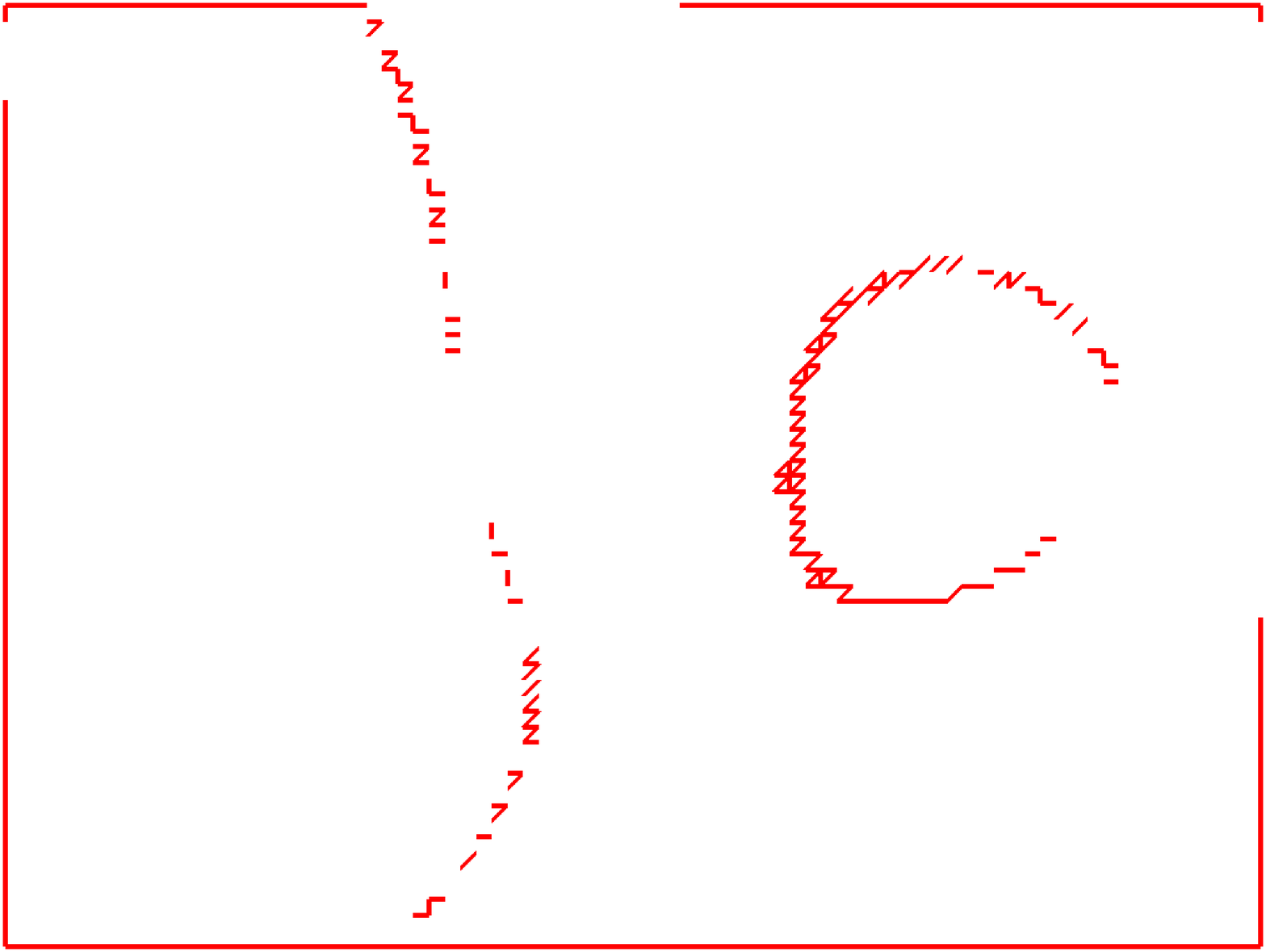}
               &\includegraphics[scale=0.08]{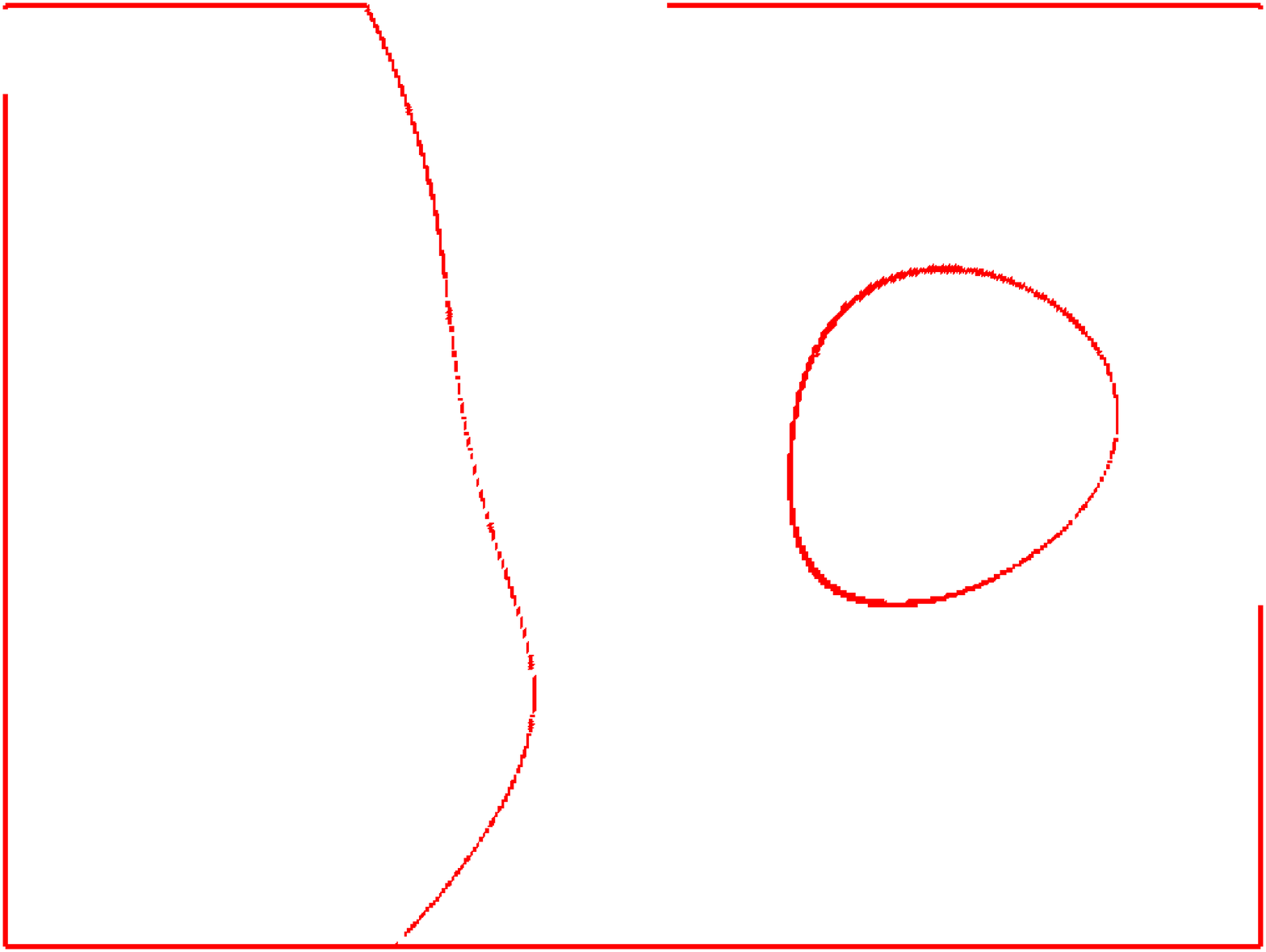}
               \\ 
 \raisebox{8mm}{$0.0001$}      
               &\includegraphics[scale=0.08]{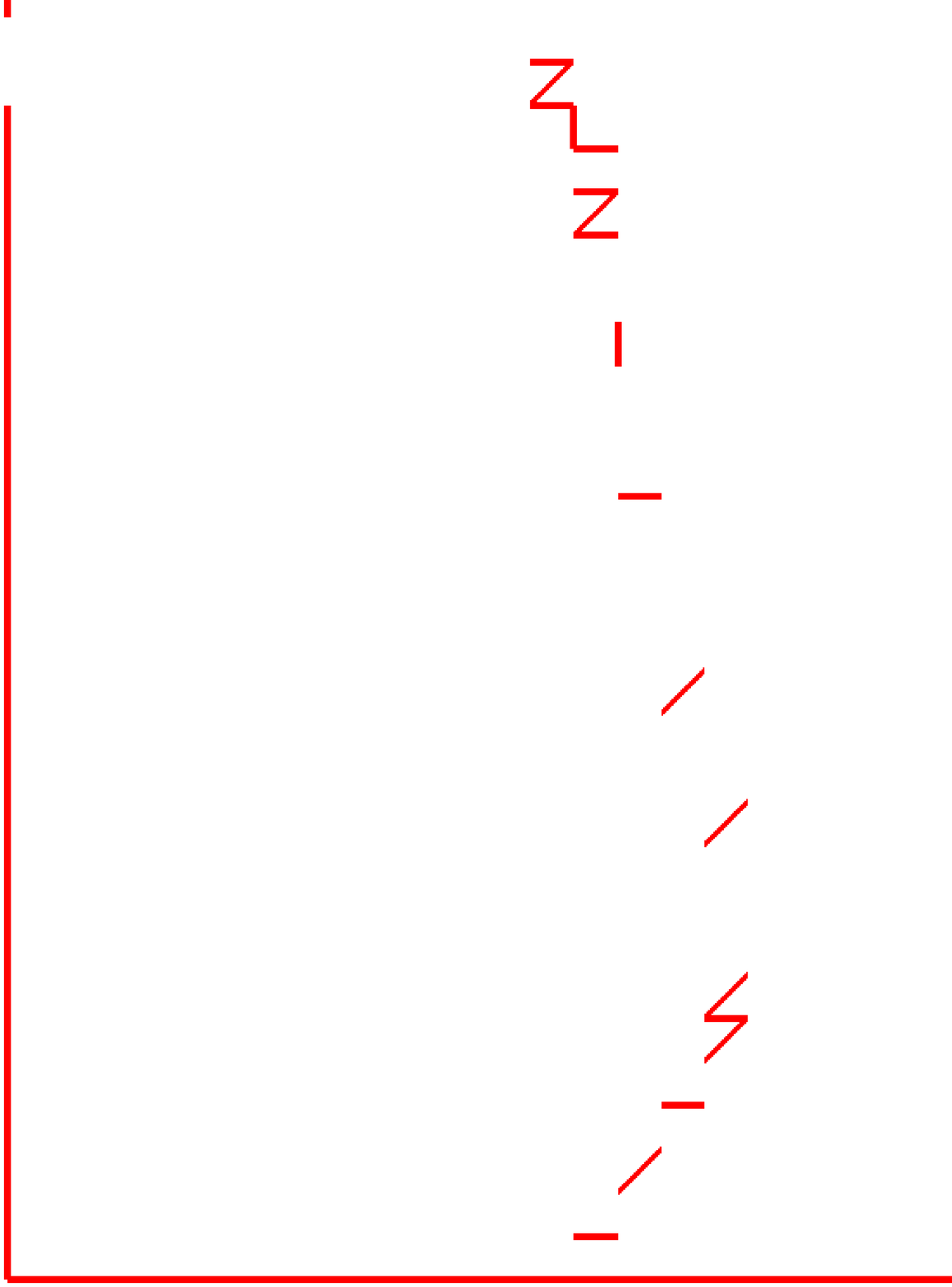}      
               &\includegraphics[scale=0.08]{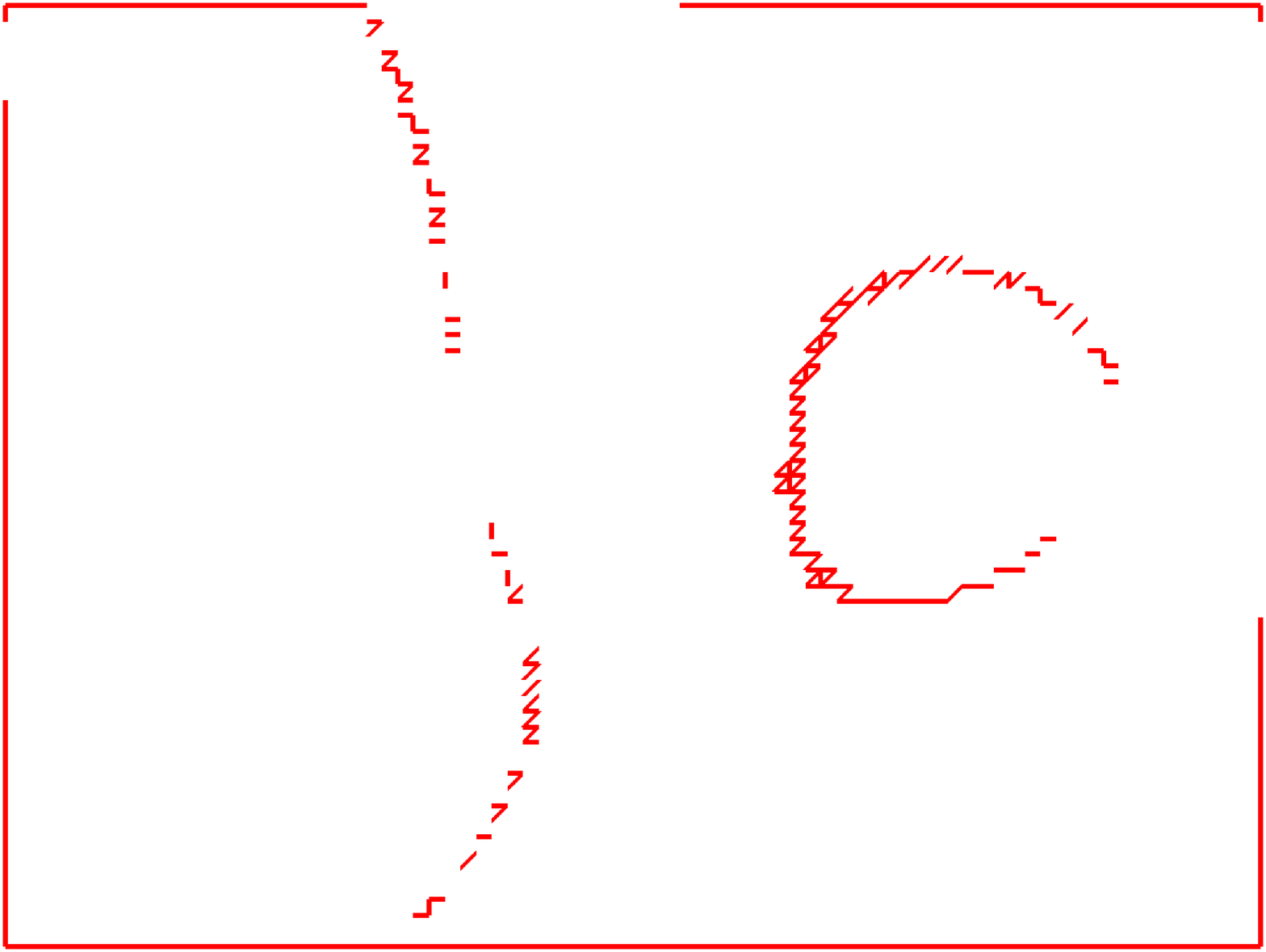}
               &\includegraphics[scale=0.08]{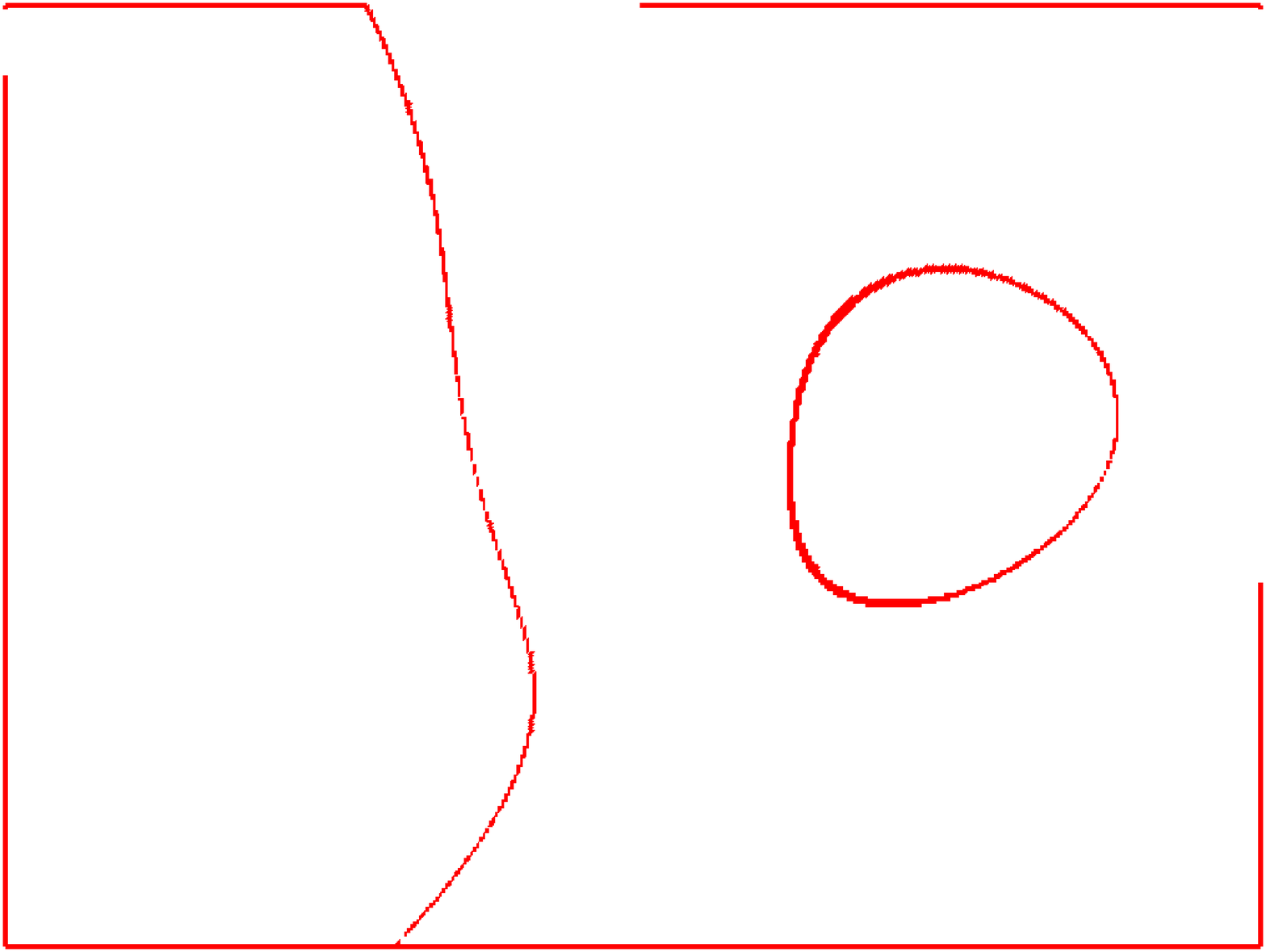}
               \\
 \raisebox{8mm}{$0.001$}       
               &\includegraphics[scale=0.08]{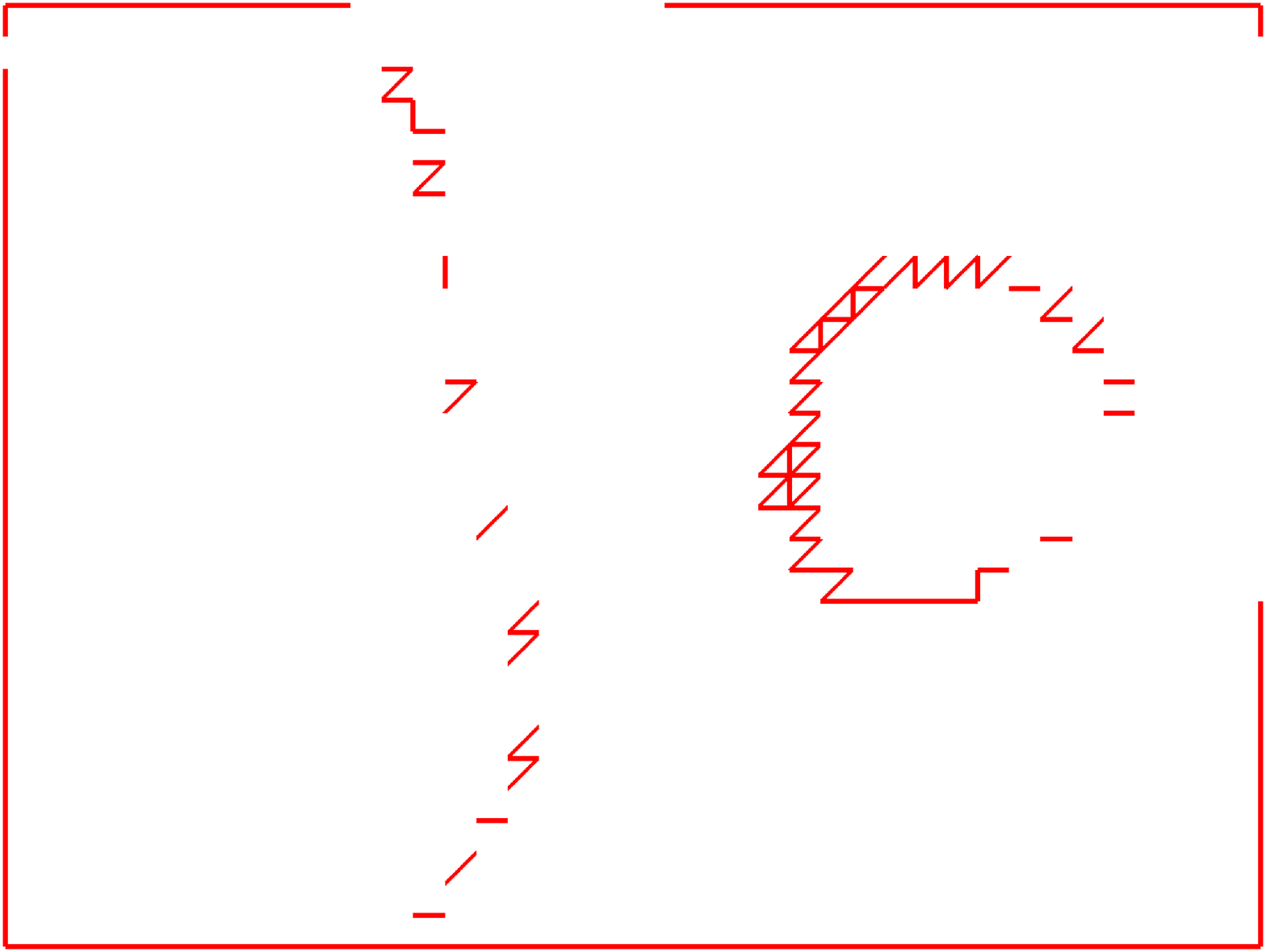}      
               &\includegraphics[scale=0.08]{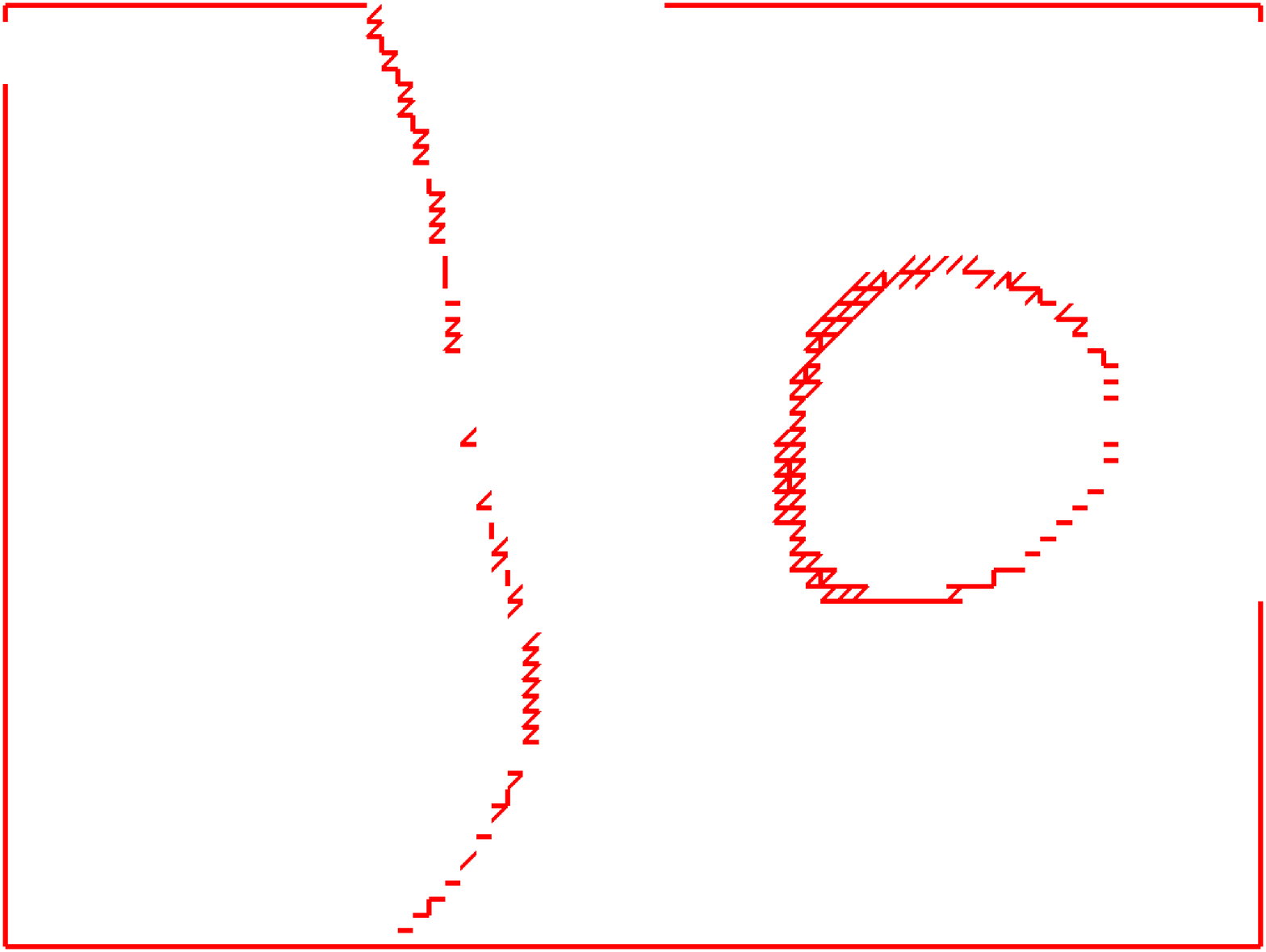}
               &\includegraphics[scale=0.08]{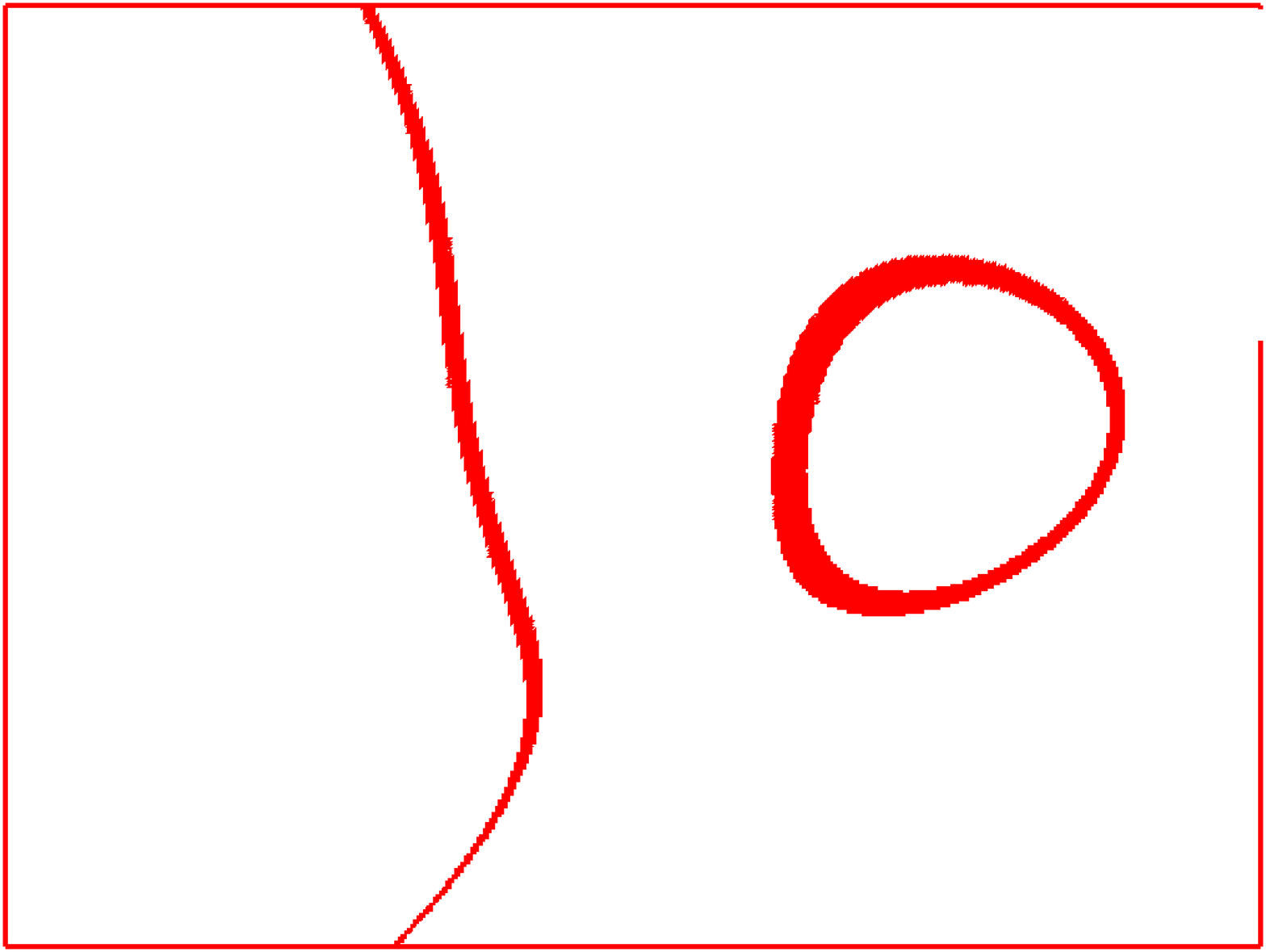}
               \\
 \raisebox{8mm}{$0.005$}       
               &\includegraphics[scale=0.08]{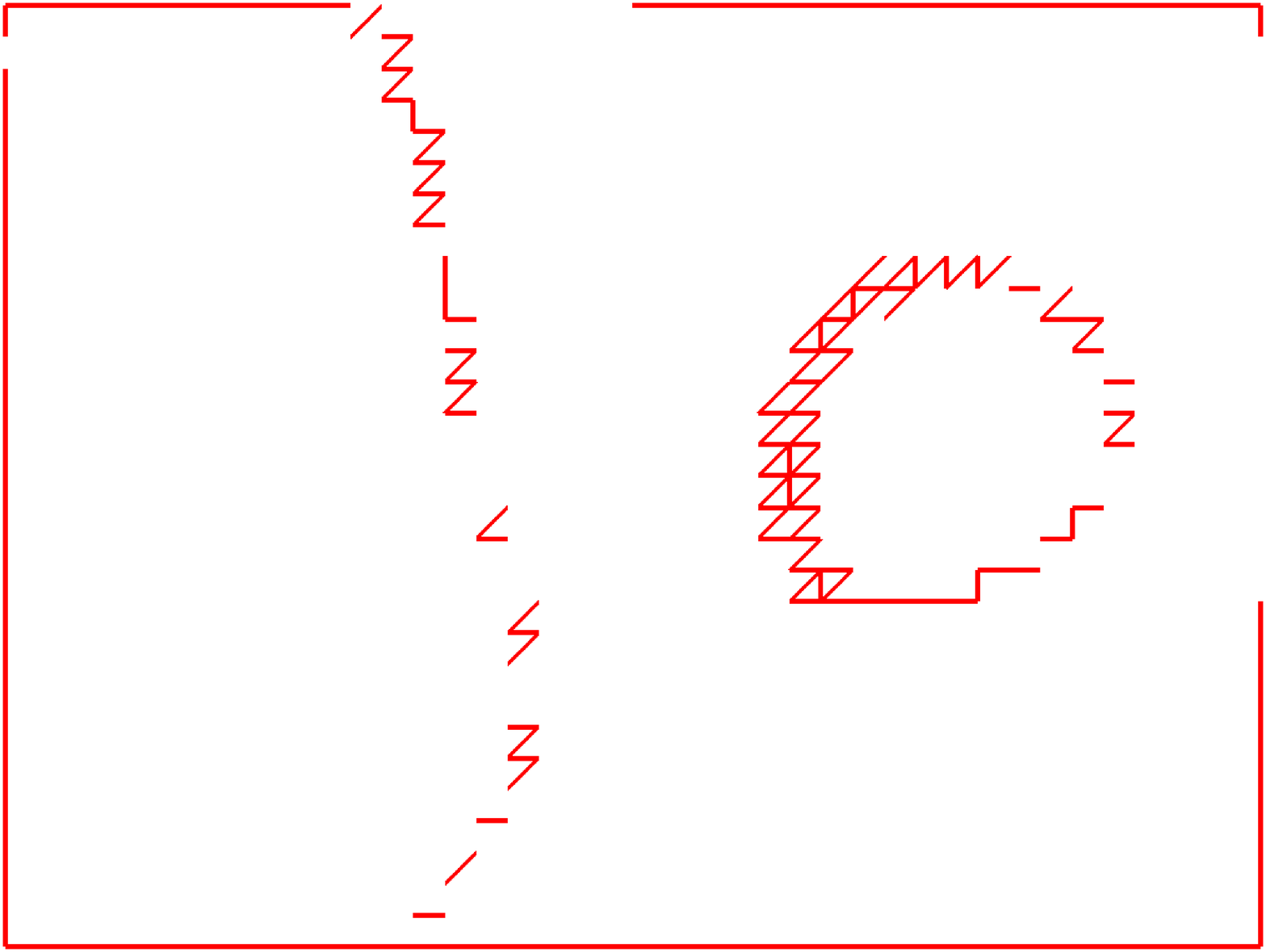}      
               &\includegraphics[scale=0.08]{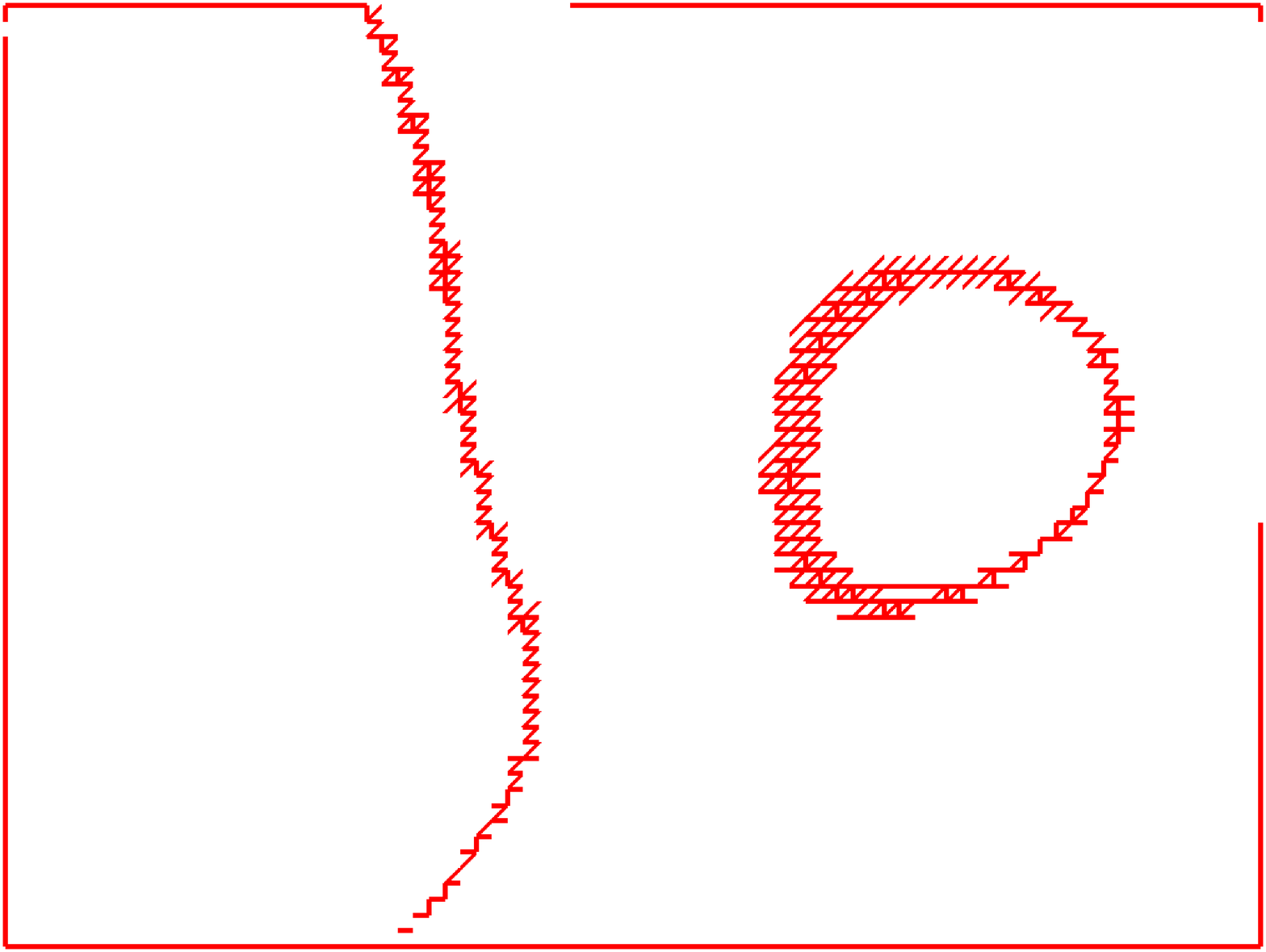}
               &\includegraphics[scale=0.08]{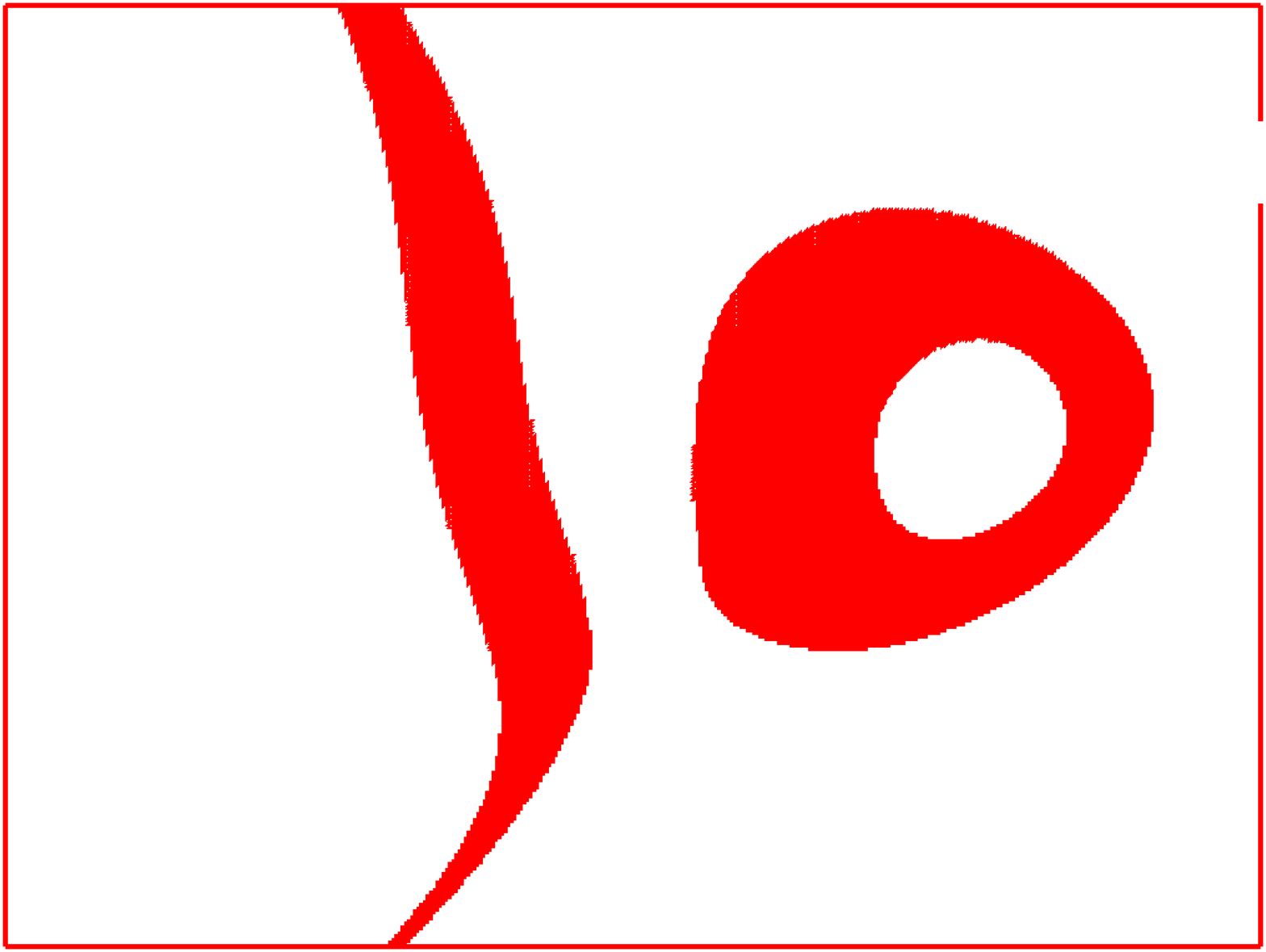}
               \\
\hline
\end{tabular}
\end{center}
\caption{Approximations for the Jacobi set in fig.~\ref{fields} for various $h$ and $\varepsilon$}
\end{table}

\newpage

\newpage
\clearpage

\begin{thebibliography}{1}

\bibitem{measure}
Edelsbrunner, H., Harer, J., Natarajan, V., Pascucci, V.: 
{Local and global comparison of continuous functions}.
\newblock In: Proceedings 16th IEEE Conference on Visualization, pp. 275--280.
IEEE Computer Society (2004).
\newblock doi:~\discretionary{}{}{}{10.1109/VISUAL.2004.68}

\bibitem{W1}
Wolpert, N.: 
{An exact and efficient approach for computing a cell in an arrangement of quadrics}.
\newblock Ph.\,D. thesis. 
Univ. des Saarlandes (2002)

\bibitem{W2}
Wolpert, N.: 
{Jacobi Curves: Computing the Exact Topology of Arrangements of Non-singular Algebraic Curves}.
\newblock In: G.~{Di Battista}, U.~Zwick (eds.) Algorithms -- ESA 2003: 11th
  Annual European Symposium, Budapest, Hungary, September 16--19, 2003.
  Proceedings, 
pp. 532--543. 
Springer, Berlin (2003).
\newblock doi:~\discretionary{}{}{}{10.1007/978-3-540-39658-1{\_}49}

\bibitem{EH}
Edelsbrunner, H., Harer, J.: 
{Jacobi Sets of multiple Morse functions}.
\newblock In: F.~Cucker, R.~DeVore, P.~Olver, E.~S{\"{u}}li (eds.) Foundations
  of Computational Mathematics, Minneapolis 2002, London Mathematical Society
  Lecture Note Series, 
pp. 37--57. 
Cambridge University Press (2004).
\newblock doi:~\discretionary{}{}{}{10.1017/CBO9781139106962.003}

\bibitem{NT}
Novikov, S., Taimanov, I.: 
{Modern Geometric Structures and Fields}.
\newblock American Mathematical Society, Providence, RI (2006)

\bibitem{Natarajan}
N, S., Natarajan, V.: 
{Simplification of Jacobi Sets}.
\newblock In: V.~Pascucci, X.~Tricoche, H.~Hagen, J.~Tierny (eds.) Topological
  Methods in Data Analysis and Visualization: Theory, Algorithms, and
  Applications, 
pp. 91--102. 
Springer, Berlin (2011).
\newblock doi:~\discretionary{}{}{}{10.1007/978-3-642-15014-2{\_}8}

\bibitem{Pascucci}
Bhatia, H., Wang, B., Norgard, G., Pascucci, V., Bremer, P.T.: 
{Local, smooth, and consistent Jacobi set simplification}.
\newblock Computational Geometry: Theory and Applications \textbf{48} (4),
311--332 (2015).
\newblock doi:~\discretionary{}{}{}{10.1016/j.comgeo.2014.10.009}

\bibitem{GFS}
{NOAA Operational Model Archive and Distribution System}. Data Transfer: NCEP GFS Forecasts (0.25 degree grid).
\newblock URL: http://nomads.ncep.noaa.gov/cgi-bin/filter\_gfs\_0p25\_1hr.pl?dir=\%2Fgfs.2018110100.

\end{thebibliography}
\end{document}